\documentclass[11pt]{article}
\usepackage{amsmath,amssymb}
\newtheorem{theo}{Theorem}[section]

\newtheorem{lemma}[theo]{Lemma}
\newtheorem{coro}[theo]{Corollary}

\newcommand{\ignore}[1]{}
\def\square{\vrule height6pt width7pt depth1pt}
\def\endpf{\hfill\square\bigskip}

\begin{document}
\title{Asymptotically optimal $K_k$-packings of dense graphs via fractional $K_k$-decompositions}
\author{Raphael Yuster
\thanks{
e-mail: raphy@research.haifa.ac.il \qquad
World Wide Web: http:$\backslash\backslash$research.haifa.ac.il$\backslash$\symbol{126}raphy}
\\ Department of Mathematics\\ University of
Haifa at Oranim\\ Tivon 36006, Israel}

\date{}

\maketitle
\setcounter{page}{1}
\begin{abstract}
Let $H$ be a fixed graph. A {\em fractional $H$-decomposition} of a graph $G$ is an assignment of nonnegative
real weights to the copies of $H$ in $G$ such that for each $e \in E(G)$, the sum of the weights of copies of $H$
containing $e$ in precisely one. An {\em $H$-packing } of a graph $G$ is a set of edge disjoint copies of $H$ in $G$.
The following results are proved. For every fixed $k > 2$, every graph with $n$ vertices and minimum degree at least
$n(1-1/9k^{10})+o(n)$ has a fractional $K_k$-decomposition and has a $K_k$-packing which covers all but $o(n^2)$
edges.
\end{abstract}

\section{Introduction}
All graphs considered here are finite, undirected and simple.
For standard graph-theoretic terminology the reader is referred to \cite{Bo}.
Let $H$ be a fixed graph. For a graph $G$, the {\em $H$-packing number}, denoted
$\nu_{H}(G)$, is the maximum number of pairwise
edge-disjoint copies of $H$ in $G$. A function $\psi$ from the set of copies of $H$ in $G$ to $[0,1]$
is a {\em fractional $H$-packing} of $G$ if $\sum_{e \in H} {\psi(H)} \leq 1$ for each $e \in E(G)$.
For a fractional $H$-packing $\psi$, let $|\psi|=\sum_{H \in {{G} \choose {H}}} \psi(H)$.
The {\em fractional $H$-packing number}, denoted $\nu^*_{H}(G)$, is defined to be the maximum value
of $|\psi|$ over all fractional packings $\psi$. Notice that, trivially, $\nu^*_{H}(G) \geq \nu_{H}(G)$.
In case $\nu_{H}(G)=e(G)/e(H)$ we say that $G$ has an {\em $H$-decomposition}.
In case $\nu^*_{H}(G)=e(G)/e(H)$ we say that $G$ has a {\em fractional $H$-decomposition}.
It is well known that computing $\nu_{H}(G)$ is NP-Hard for every fixed graph $H$ with more than two edges in some connected
component \cite{DoTa}.
It is well known that computing $\nu^*_{H}(G)$ is solvable in polynomial time for every fixed graph $H$ as this amounts to solving a
(polynomial size) linear program.

The combinatorial aspects of the $H$-packing and $H$-decomposition problems have been studied extensively.
Wilson in \cite{Wi} has proved that whenever $n \geq n_0=n_0(H)$, and $K_n$ satisfies two obvious necessary divisibility
requirements, $K_n$ has an $H$-decomposition. The $H$-packing problem for $K_n$ ($n \geq n(H)$) was solved \cite{CaYu},
by giving a closed formula for $\nu(H,K_n)$. For graphs $G$ which are not complete, giving sufficient conditions
guaranteeing an $H$-decomposition or, at least, a packing covering all but a small fraction of the edges seems to be an
extremely difficult task, even if we assume that $G$ is dense. To be more precise, let us define the following problem and parameter.

{\bf Problem:}\, Given a fixed graph $H$, determine $c_H$, which is the supremum of all possible
constants $c$ guaranteeing
that every graph $G$ with $n$ vertices and minimum degree $\delta(G) \geq (1-c)(n-1)$ has $\nu_H(G) \geq (1-o_n(1))e(G)/e(H)$.

Wilson's Theorem implies that $c_H \geq 0$ exists. In fact, $c_H \geq 0$ can also be derived from R\"odl's result \cite{Ro}
which is weaker than Wilson's for graphs, but is more general since it applies to $r$-uniform hypergraphs as well.
Gustavsson \cite{Gu} has proved that $c_H > 0$ for every graph $H$. However, Gustavsson's lower bound for $c_H$ is
horribly close to $0$. Already for $H=K_3$ it only gives $c_{K_3} > 10^{-24}$, and, more generally, if $H$ has
$k$ vertices then $c_H > 10^{-37}k^{-94}$. Gustavsson's result does however, show that the minimum degree
requirement, together with necessary divisibility conditions, guarantees an $H$-decomposition. The exact value of $c_H$ is
unknown for any fixed non-bipartite graph $H$. Notice that, trivially, $c_H=1$ in case $H$ is a bipartite graph as
the Tur\'an number of such graphs is $o(n^2)$. However, if we insist on having a packing of size $\lfloor e(G)/e(H) \rfloor$
then it is known that a minimum degree of $0.5n(1+o_n(1))$ suffices for each fixed bipartite graph $H$ (other than the trivial $K_2$)
having a vertex of degree one (this includes all trees) and this is asymptotically tight \cite{Yu}.

In this paper we prove the first reasonable general lower bound for $c_H$.
\begin{theo}
\label{t1}
For all $k \geq 3$,  $c_{K_k} \geq 1/9k^{10}$.
\end{theo}
Although Theorem \ref{t1} is stated only for $K_k$, a simple argument given in the sequel shows that the same
lower bound holds for any graph $H$ with $k$ vertices.
Theorem \ref{t1} is deduced as a combination of two powerful theorems, the first of which is the following.
\begin{theo}
\label{t2}
For all $k \geq 3$,  any graph $G$ with $n$ vertices and $\delta(G) \geq n(1-1/9k^{10})+o(n)$ has a fractional $K_k$-decomposition.
\end{theo}
The proof of Theorem \ref{t2} appears in the following section.

Recently, Haxell and R\"odl \cite{HaRo} proved that the $H$-packing number and the fractional $H$-packing number are very close
for dense graphs. A simpler and more general proof of their result appears in \cite{Yu2}.
\begin{theo}
\label{t3}[Haxell and R\"odl \cite{HaRo}]
For any fixed graph $H$, if $G$ has $n$ vertices then $\nu^*_{H}(G)-\nu_{H}(G) = o(n^2)$. \endpf
\end{theo}

Theorem \ref{t2} gives that for sufficiently large $n$, any graph $G$ with $n$ vertices and $\delta(G) \geq n(1-1/9k^{10})+o(n)$
has $\nu^*_{K_k}(G)=e(G)/e(K_k)$. Thus, by Theorem \ref{t3}, it also has
$\nu_{K_k}(G) \geq e(G)/e(K_k)-o(n^2)=(1-o_n(1))e(G)/e(K_k)$. Consequently, $c_{K_k} \geq 1/9k^{10}$ and Theorem \ref{t1}
follows.

Finally, we note that an $1/(k+1)$ upper bound for $c_{K_k}$ is given in the final section together with some additional
concluding remarks.

\section{Proof of Theorem \ref{t2}}
Let ${\cal F}$ be a fixed family of graphs. An {\em ${\cal F}$-decomposition} of a graph $G$ is a set $L$ of subgraphs
of $G$, each isomorphic to an element of ${\cal F}$, and such that each edge of $G$ appears in precisely one element of $L$.
Let $K_t^-$ denote the complete graph with $t$ vertices, missing one edge.
Let ${\cal F}_k=\{K_k~,~ K_{2k-1}~,~ K_{2k-1}^-\}$.
The proof of Theorem \ref{t2} is a corollary of the following stronger theorem.
\begin{theo}
\label{t21}
For all $k \geq 3$, every graph with $n$ vertices and minimum degree
at least $n(1-1/9k^{10})+o(n)$ has an ${\cal F}_k$-decomposition.
\end{theo}
The following simple lemma shows that Theorem \ref{t2} is a corollary of Theorem \ref{t21}.
\begin{lemma}
\label{l22}
For all $k \geq 2$, the graphs $K_{2k-1}$ and  $K_{2k-1}^-$ have a fractional $K_k$-decomposition.
\end{lemma}
{\bf Proof:}\,
It is trivial that for all $k' \geq k$, $K_{k'}$ has a fractional $K_k$-decomposition. In Particular,
$K_{2k-1}$ has a fractional $K_k$-decomposition.
Let $A=\{u,v\}$ denote the set of the two non-adjacent vertices of $K_{2k-1}^-$,
and let $B$ denote the set of the remaining $2k-3$ vertices. Each edge incident with $A$ lies on
${{2k-4} \choose {k-2}}$ copies of $K_k$. Each edge with both endpoints in $B$ lies
on $2{{2k-5} \choose {k-3}}$ copies of $K_k$ that contain a vertex of $A$. Since
${{2k-4} \choose {k-2}}=2{{2k-5} \choose {k-3}}$, by assigning the value
$1/{{2k-4} \choose {k-2}}$ to each copy of $K_k$ containing a vertex of $A$, and assigning
the value $0$ to the remaining copies of $K_k$, we obtain a fractional $K_k$-decomposition of $K_{2k-1}^-$.
\endpf

We now focus on proving Theorem \ref{t21}. For the rest of this section, let $t=2k-1$.
Notice that the $o(n)$ term in the statement of Theorem \ref{t21} allows us to assume, whenever necessary,  that $n$ is sufficiently
large.
In the proof of Theorem \ref{t21} it will be convenient to use Wilson's Theorem \cite{Wi} mentioned in the introduction.
(We note that it is also possible to use R\"odl's packing theorem \cite{Ro} instead of Wilson's Theorem at the price of some
complication in the proof). Wilson's Theorem applied to $K_t$ states the following.
\begin{lemma}
\label{l23}[Wilson]
Let $t > 2$ be a positive integer. There exists $N=N(t)$ such that for all $n > N$ with $n \equiv 1,t \bmod t(t-1)$,
there is a decomposition of $K_n$ into $K_t$.
\end{lemma}

We shall prove Theorem \ref{t21} under the relaxed assumption that $n \equiv 1 \bmod t(t-1)$.
We first need to justify this relaxation.
Indeed, if $1 < b \leq t(t-1)$ and $n \equiv b \bmod t(t-1)$ then we can perform the following preprocessing.
Let $v$ be any vertex of $G$, and let $N(v)$ denote its neighborhood. Let $G[N(v)]$ be the subgraph induced by
this neighborhood. Notice that $G[N(v)]$ has less than $n$ vertices, but has minimum degree at least
$n(1-2/9k^{10})+o(n)$. By the theorem of Hajnal and Szemer\'edi \cite{HaSz}, such a high minimum degree for a graph with at
most $n$ vertices is far more than what is needed in order to guarantee that $G[N(v)]$ has a spanning subgraph $G'$
with each component of $G'$ being either a $K_{k-1}$ or a $K_{3k-4}$.
In fact, a minimum degree of at least $n(1-1/(3k-4))$ already guarantees the existence of such a $G'$.
Now, if $H$ is a $K_{k-1}$ component of $G'$ then $H \cup \{v\}$ is a $K_k$ copy of $G$. If $H$ is a $K_{3k-4}$
component of $G'$ then $H \cup {v}$ contains two edge-disjoint subgraphs, one being $K_k$ and the other being
$K_t^-$, and with all $3k-4$ edges between $v$ and the vertices of $H$ absorbed.
We thus have a set of edge-disjoint subgraphs of $G$, each being either $K_k$ or $K_t^-$, and that absorb all
edges incident with $v$. Deleting $v$ and the edges of these subgraphs we remain with a graph with $n-1$ vertices
and minimum degree at least $\delta(G) - (3k-4) \geq (n-1)(1-1/9k^{10})+o(n-1)$. Repeating this process at most
$b-1$ times we eventually have a graph with $n'=n-b+1 \equiv 1 \bmod t(t-1)$ vertices and minimum degree
at least $\delta(G)- (3k-4)((2k-1)(2k-2)-1) \geq n'(1-1/9k^{10})+o(n')$, which satisfied our relaxed assumption.
Our preprocessing shows that any ${\cal F}_k$-decomposition of this resulting $n'$-vertex graph can be extended to
an ${\cal F}_k$-decomposition of the original $n$-vertex graph.

\bigskip
\noindent
{\bf Proof of Theorem \ref{t21}}
We may assume that $n \equiv 1 \bmod t(t-1)$ and that, whenever necessary, $n$ is sufficiently large as a function of $k$ (and hence $t$).
In particular, $n > N(t)$ where $N(t)$ is the constant from Lemma \ref{l23}.
Fix a $K_t$-decomposition of $K_n$. Namely, let $D$ be a
family of $t$-sets of $[n]=\{1,\ldots,n\}$ such that each pair appears in precisely one element of $D$. Such a $D$ is also
called a $t$-design. Notice that $|D|={n \choose 2}/{t \choose 2}$. For a permutation $\pi$ of $[n]$, and
for $S \in D$, let $S_\pi = \{\pi(j) ~: ~ j \in S\}$. Hence, $D_\pi = \{S_\pi ~:~ S \in D\}$ is also a $t$-design.
Let $G$ be an $n$-vertex graph with vertex set $[n]$.
For a permutation $\pi$ of $[n]$ let $G_\pi$ be the family of  ${n \choose 2} / {t \choose 2}$ edge-disjoint subgraphs
of $G$ whose elements are the induced subgraphs of $G$ on $S_\pi$, for all $S \in D$. Notice that if $G=K_n$ then,
trivially, $G_\pi$ is a $K_t$-decomposition for each $\pi$, but if $G \neq K_n$, $G_\pi$ contains elements that are
not isomorphic to $K_t$. The following is a simple corollary of Lemma \ref{l23}.
\begin{coro}
\label{c24}
Let $0 < \alpha < 1$ be fixed.
Let $G$ be a graph with $n > N(t)$ vertices, $n \equiv 1 \bmod t(t-1)$. If $\delta(G) \geq (1-\alpha)n$
and $\pi$ is any permutation of $[n]$ then $G_\pi$ has at least $(1-o(1))n^2(\frac{1}{t(t-1)}-\frac{\alpha}{2})$ elements
isomorphic to $K_t$ and at most $(1+o(1))n^2(\frac{\alpha}{2})({t \choose 2}-1)$ edges appear in elements of $G_\pi$ that are
not isomorphic to $K_t$.
\end{coro}
{\bf Proof:}\,
The number of non-edges of $G$ is at most ${n \choose 2}-(1-\alpha) n^2/2$. Thus, $G_\pi$ has at most
${n \choose 2}-(1-\alpha) n^2/2$ elements that are not $K_t$ and therefore $G_\pi$ has at least
${n \choose 2}/{t \choose 2} - {n \choose 2} + (1-\alpha) n^2/2$ elements isomorphic to $K_t$ and at most
$({n \choose 2}-(1-\alpha) n^2/2)({t \choose 2}-1)$ edges of $G$ are in non-$K_t$ elements of $G_\pi$. \endpf

Assume that $\delta(G) \geq n(1-1/9k^{10}) + o(n)$.
Our goal is to show that there exists a permutation $\pi$, such that $G_\pi$ has some ``nice'' properties.
Let $A_\pi$ denote the set of edges of $G$ that appear in non-$K_t$ elements of $G_\pi$.
By Corollary \ref{c24}, with $\alpha=1/9k^{10}$,
\begin{equation}
\label{e0}
|A_\pi| \leq (1+o(1))n^2\left(\frac{1}{18k^{10}}\right)\left({t \choose 2}-1\right) \leq (1+o(1))\frac{n^2}{9k^8}.
\end{equation}
Consider the spanning subgraph of $G$ consisting of the edges of $A_\pi$. It will not be confusing to denote this subgraph by
$A_\pi$ as well. Let $F_\pi \subset G_\pi$ be the set of $K_t$-elements of $G_\pi$.
Put $r={k \choose 2}-1$.
We say that an $r$-subset $S=\{H_1,\ldots,H_r\}$ of $F_\pi$ is {\em good for $e \in A_\pi$} if we can select edges
$f_i \in H_i$ such that $\{f_1,\ldots,f_r,e\}$ is the set of edges of a $K_k$ in $G$.
We say that $\pi$ is {\em good} if for each $e \in A_\pi$ there exists an $r$-subset $S(e)$ of $F_\pi$ such that
$S(e)$ is good for $e$ and such that if $e \neq e'$ then $S(e) \cap S(e') = \emptyset$.

\begin{lemma}
\label{l25}
If $\pi$ is good then $G$ has an ${\cal F}_k$-decomposition.
\end{lemma}
{\bf Proof:}\,
For each $e \in A_\pi$, pick a copy of $K_k$ in $G$ containing $e$ and precisely one edge from each element of $S(e)$.
As each element of $S(e)$ is a $K_t$, deleting one edge from such an element results in a $K_t^-$.
We therefore have $|A_\pi|$ copies of $K_k$ and $|A_\pi|({k \choose 2}-1)$
copies of $K_t^-$, all being edge disjoint. The remaining element of $F_\pi$ not belonging to
any of the $S(e)$ are each a $K_t$, and they are edge-disjoint from each other and from the previously selected $K_k$ and $K_t^-$.
\endpf

Our goal in the remainder of this section is, therefore, to show that there exists a good $\pi$.
We use probabilistic and counting arguments to derive this fact. We will show that with positive probability,
a randomly selected $\pi$ is good. We begin by showing that with high probability, a randomly
selected $\pi$ has the property that $A_\pi$ has a relatively small maximum degree.

Let $v$ be any vertex of $G$, and let $E(v)$ denote the set of edges incident with $v$.
By our assumption, $|E(v)| \geq (1-1/9k^{10})n$.
Let $E_\pi(v) = E(v) \cap A_\pi$. Notice that if $\pi$ is selected uniformly at random then $|E_\pi(v)|$, which is the degree of $v$ in
$A_\pi$, is a random variable.
We say that a subset $S \subset E(v)$ is {\em separated by $\pi$} if each edge of $S$ belongs to a different element of
$G_\pi$. Let $\beta=3/k^8$ and consider any fixed set $S \subset E(v)$ with $|S|=\lfloor \beta n \rfloor$.
We shall prove that the probability that $S \subset E_\pi(v)$ {\em and} that $S$ is separated by $\pi$ is much smaller 
than the total number of subsets of $E(v)$ with size $\lfloor \beta n \rfloor$.
Thus, the probability that the degree of $v$ in $A_\pi$ exceeds $t \beta n$ is also very small (in fact, much smaller than $1/n$),
and consequently, the maximum degree of $A_\pi$ is at most $t\beta n$ almost surely.
Let $S=\{e_1,\ldots, e_m\}$ where $m=\lfloor \beta n \rfloor$ and let $e_i=(v,v_i)$.
Let $S_\pi(i)$ denote the element of $G_\pi$ to which the edge $e_i$ belongs, $i=1,\ldots,m$.
Notice that $S$ is separated by $\pi$ if and only if $S_\pi(i) \neq S_\pi(j)$ for all $i \neq j$.
Clearly, by using conditional probabilities we have
\begin{equation}
\label{e1}
\Pr[(S \subset E_\pi(v)) \wedge (S {\rm~is~separated~by~\pi})] =
\end{equation}
$$
\prod_{i=1}^m \Pr[(e_i \in A_\pi) \wedge (\forall j < i,~ S_\pi(i) \neq S_\pi(j))  ~|~
(\{e_1,\ldots,e_{i-1}\} \subset A_\pi) \wedge (S_\pi(j) \neq S_\pi(j'),~ 1 \leq j <j' <i)].
$$
We shall prove that each term in the product appearing in the r.h.s. of the last equation is small.
For this purpose we require a lemma which quantifies the fact that in a graph with high minimum degree
every edge appears on many copies of $K_t$.
\begin{lemma}
\label{l26}
If $G^*$ is a graph with $n^*$ vertices and minimum degree at least $n^*-r$ then every edge of $G^*$ appears
on at least $\frac{1}{(t-2)!}\prod_{i=2}^{t-1}(n^*-ir)$ distinct copies of $K_t$.
\end{lemma}
We prove the lemma by induction on $t$. For $t=2$ and $t=3$ the lemma is obvious. Assume the lemma holds for
all $t' < t$. Let $e=(u,v)$ be an edge of $G^*$.
Let $N(u,v)$ denote the set of common neighbors of $u$ and $v$. Clearly, $|N(u,v)| \geq n^*-2r$.
Let $G^{**}=G^*[N(u,v)]$. The minimum degree of $G^{**}$ is at least $|N(u,v)|-r$.
It follows that $G^{**}$ has at least $(n^*-2r)(n^*-3r)/2$ edges. The number of distinct copies of $K_{t-2}$ in $G^{**}$
is equal to the number of distinct copies of $K_t$ containing $e$ in $G^*$.
By the induction hypothesis, each
edge of $G^{**}$ appears in at least $\frac{1}{(t-4)!}\prod_{i=4}^{t-1}(n^*-ir)$ distinct copies of $K_{t-2}$.
Since each copy of $K_{t-2}$ is counted $(t-2)(t-3)/2$ times the number of distinct
copies of $K_{t-2}$ in $G^{**}$ is at least
$$
\frac{(n^*-2r)(n^*-3r)}{2} \frac{1}{(t-4)!}\prod_{i=4}^{t-1}(n^*-ir) \frac{2}{(t-2)(t-3)}=
\frac{1}{(t-2)!}\prod_{i=2}^{t-1}(n^*-ir)
$$
as required. \endpf

\begin{coro}
\label{c27}
If $G^*$ is a graph with $n^*$ vertices and minimum degree at least $n^*(1-\gamma)$ then,
for every edge $e$ of $G^*$, the probability that a randomly selected $t$-vertex subgraph of $G^*$ that contains $e$
is not a $K_t$ is at most $1-(1-t\gamma)^{t-2}$.
\end{coro}
{\bf Proof:}\,
There are precisely ${{n^*-2} \choose {t-2}}$ subgraphs with $t$ vertices that contain the edge $e$.
By Lemma \ref{l26}, with $r=\gamma n^*$, the number of $K_t$-subgraphs that contain $e$ is at least
$$
\frac{1}{(t-2)!}\prod_{i=2}^{t-1}(n^*-i\gamma n^*) >
\frac{1}{(t-2)!}(n^*)^{t-2}(1-t\gamma)^{t-2} > {{n^*-2} \choose {t-2}}(1-t\gamma)^{t-2}.
$$
Thus, the probability that a randomly selected $t$-vertex subgraph of $G^*$ that contains $e$
is not a $K_t$ is at most $1-(1-t\gamma)^{t-2}$. \endpf

Corollary \ref{c27} enables us to estimate the terms in the r.h.s. of (\ref{e1}).
Let $Y_\pi(j)$ be the set of $t$ vertices of the element $S_\pi(j)$. Notice that given the knowledge that
$S_\pi(j) \neq S_\pi(j')$ for $1 \leq j < j' < i$ implies, in particular, the knowledge that $Y_\pi(j) \cap Y_\pi(j')=\{v\}$.
Thus, if $W = \cup_{j=1}^{i-1} Y_\pi(j)$ then $|W|=(i-1)(t-1)+1$. Thus, we know the size of $W$.
In order to prove an upper bound on each term of the r.h.s. of (\ref{e1}) it suffices to prove an upper bound
on $\Pr[(e_i \in A_\pi)  \wedge (S_\pi(i) \cap W = \{v\}) | W] $ whose value does not depend on the specific set $W$, but which may,
and will, depend on the prior knowledge that $|W|=(i-1)(t-1)+1$.
Indeed, let $G^*$ be the induced subgraph of $G$ obtained by deleting the set of vertices $W - \{v\}$. Notice that $G^*$ has
$n^*=n-(i-1)(t-1)$ vertices. Clearly, $\Pr[(e_i \in A_\pi)  \wedge (S_\pi(i) \cap W = \{v\}) | W] $ is precisely the
probability that a randomly selected $t$-vertex subgraph of $G^*$ containing $e_i$ is not a $K_t$.
Using the fact that $(i-1)(t-1) < mt << n/2$, we have that the minimum degree of $G^*$ is a least $n^*-n/9k^{10} \geq n^*(1-2/9k^{10})$.
Using $\gamma=2/9k^{10}$, we have by Corollary \ref{c27} that
$\Pr[(e_i \in A_\pi)  \wedge (S_\pi(i) \cap W = \{v\}) | W]  \leq 1-(1-t\gamma)^{t-2}$.
Consequently each term in (\ref{e1}) is bounded from above by $1-(1-t\gamma)^{t-2}$.
We therefore have
\begin{equation}
\label{e2}
\Pr[(S \subset E_\pi(v)) \wedge (S {\rm~is~separated~by~\pi})] \leq (1-(1-t\gamma)^{t-2})^m.
\end{equation}

\begin{lemma}
\label{l28}
With probability $1-o(1)$, $A_\pi$ has maximum degree at most $6n/k^{7}$.
\end{lemma}
{\bf Proof:}\,
Since the maximum degree of each element of $G_\pi$ is at most $t-1$ we have, by the definition of $E_\pi(v)$, that there
is $S \subset E_\pi(v)$ such that $S$ is separated by $\pi$ and $|S| \geq |E_\pi(v)|/(t-1)$.
Thus, it suffices to show that for each vertex $v$, $E_\pi(v)$ has no subset separated by $\pi$ of size greater than $3n/k^8$
with probability $1-o(1/n)$ since this implies that $|E_\pi(v)| \leq (t-1)3n/k^8 \leq 6n/k^7$ with probability
$1-o(1/n)$ and hence the maximum degree of $A_\pi$ is at most $6n/k^7$ with probability $1-o(1)$.
Indeed, by (\ref{e2}) the probability that $E_\pi(v)$ has a subset separated by $\pi$ of size $m=\lfloor \beta n \rfloor=
\lfloor 3n/k^8 \rfloor$ is at most
${{n-1} \choose m}(1-(1-t\gamma)^{t-2})^m$ where $\gamma=2/9k^{10}$.
We therefore have
$$
{{n-1} \choose m}(1-(1-t\gamma)^{t-2})^m \leq
\left(\frac{1}{\beta^\beta (1-\beta)^{1-\beta}} \left(1-(1-\frac{4}{9k^9})^{2k-3}\right)^\beta\right)^n
$$
$$
\leq
\left(\frac{1}{\beta^\beta (1-\beta)^{1-\beta}} \left(1-(1-\frac{1}{k^8})\right)^\beta\right)^n=
\left(\frac{1}{3^\beta (1-\beta)^{1-\beta}} \right)^n=o\left(\frac{1}{n}\right).
$$
\endpf

By Lemma \ref{l28}, we may {\em fix} a permutation $\pi$ for which $A_\pi$ has maximum degree at most
$6n/k^7$. Let $A_\pi=\{e_1,\ldots,e_m\}$. We perform the following algorithm which has $m$ iterations.
In the $i'$th iteration we pick an $r$-subset $S(e_i)$ of $F_\pi$ which is good for $e_i$
and which satisfies the following two properties:
\begin{enumerate}
\item
$S(e_i) \cap S(e_j) = \emptyset$ for all $j=1,\ldots,i-1$.
\item
For each $v$, let $f_i(v)$ denote the number of edges incident with $v$ and which belong to some element of $S(e_j)$,
$j \leq i$, and where $v$ is not an endpoint of $e_j$. Then, $f_i(v) \leq n/(2k)$.
\end{enumerate}
Notice that if we can complete all $m$ iterations of the algorithm then the first requirement guarantees that
$\pi$ is good and hence, by Lemma \ref{l25} we are done. The second requirement is needed in order to guarantee that
the algorithm will, indeed, complete all $m$ iterations. We therefore need to prove the following lemma.
\begin{lemma}
\label{l29}
If $A_\pi$ has maximum degree at most $6n/k^7$ then the algorithm completes all $m$ iterations.
\end{lemma}
{\bf Proof:}\,
The most difficult case is to prove that the $m$'th iteration can also be completed, assuming all previous iterations
have completed. Let $e_m=(u,v)$. We define several parameters.
Let $a_1$ denote the number of $K_k$ copies of $G$ that contain $e_m$.
Let $a_2$ denote the number of $K_k$ copies of $G$ that contain $e_m$ and also contain two edges from the same
element of $F_\pi$. Let $a_3$ denote the number of $K_k$ copies of $G$ that contain $e_m$ and also contain another
edge from $A_\pi$. Let $E_m$ denote the set of $(m-1)r{t \choose 2}$ edges in all elements of $\cup_{i=1}^{m-1}S(e_i)$.
Let $a_4$ denote the number of $K_k$ copies of $G$ that contain $e_m$ and also contain an edge of $E_m$.
Let $V_m$ be the subset of vertices of $G$ where $x \in V_m$ if and only if  $x \neq u,v$ and $f_{m-1}(x) > n/k-r(t-1)$.
Let $F_m$ be the set of all edges of all copies of $F_\pi$ that contain at least one vertex of $V_m$.
Let $a_5$ denote the number of $K_k$ copies of $G$ that contain $e_m$ and also contain an edge of $F_m$.

We claim that if $a_1 > a_2+a_3+a_4+a_5$ then the $m$'th iteration can be completed. Indeed, if this is the case then
by the definitions of $a_2$ and $a_3$ there exists a copy of $K_k$ in $G$ which contain $e_m$, and whose other $r$
edges all belong to distinct elements of $F_\pi$, say, $S(e_m)=\{H_1,\ldots,H_r\}$. Furthermore, by the definition of $a_4$ we may
also assume that no $H_i$ is an element of a previous $S(e_j)$ for $j < m$, and hence $S(e_m) \cap S(e_j) = \emptyset$
for all $j=1,\ldots,m-1$. Finally, by the definition of $a_5$ we may assume that no $H_i$ contains a vertex of
$V_m$. Thus, for $x \in V_m$ we have $f_m(x)=f_{m-1}(x) \leq n/(2k)$ by our assumption.
By definition, since $u,v$ are incident with $e_m$ we have $f_m(v)=f_{m-1}(v) \leq n/(2k)$ and $f_m(u)=f_{m-1}(u) \leq n/(2k)$.
For $x \notin V_m \cup \{u,v\}$ notice that we have $f_m(x) \leq f_{m-1}(x)+r(t-1) \leq n/(2k)$ as well.

It remains to show that $a_1 > a_2+a_3+a_4+a_5$. We now estimate these parameters.
A similar proof to that of Corollary \ref{c27} (where we use $n$ instead of $n^*$, $k$ instead of $t$ and $\gamma=1/9k^{10}$
immediately gives
\begin{equation}
\label{e3}
a_1 \geq {{n-2} \choose {k-2}}(1-\frac{1}{9k^9})^{k-2} \geq n^{k-2}\frac{0.95}{(k-2)!}(1-o(1)).
\end{equation}
Consider a pair of edges $f_1,f_2$ that belong to the same element of $F_\pi$. If they are both independent from $e_m$ then
$\{e_m,f_1,f_2\}$ spans at least five vertices. Thus, there are at most ${{n-5} \choose {k-5}}$ copies of $K_k$ that
contain all three of them. As there are less than $|F_\pi|t^4=\Theta(n^2)$ possible choices for pairs $f_1,f_2$
the overall number of such copies is $O(n^{k-3})$. If $f_1,f_2$ are not independent from $e_m$ then
we must have that $\{e_m,f_1,f_2\}$ spans at least four vertices. Thus, there are at most ${{n-4} \choose {k-4}}$ copies of
$K_k$ that contain all three of them. However, there are less than $2n$ edges not independent from $e_m$,
so the overall number of choices for $f_1,f_2$ is only $O(n)$. Overall there are, again, only $O(n^{k-3})$ such copies.
We have proved that
\begin{equation}
\label{e4}
a_2 = O(n^{k-3}).
\end{equation}
Consider an edge $f \in A_\pi$ with $f \neq e_m$. If $f$ and $e_m$ are independent then there are at most
${{n-4} \choose {k-4}}$ copies of $K_k$ containing both of them. Overall, there are at most
$m{{n-4} \choose {k-4}}$ such copies. If $f$ and $e_m$ are not independent then there are at most
${{n-3} \choose {k-3}}$ copies of $K_k$ containing both of them. However, the maximum degree of $A_\pi$ is
at most $6n/k^7$ and hence there are at most $12n/k^7$ choices for $f$. Using (\ref{e0}) we therefore have
\begin{equation}
\label{e5}
a_3 \leq m{{n-4} \choose {k-4}} + \frac{12n}{k^7}{{n-3} \choose {k-3}}
\leq (1+o(1))\frac{n^2}{9k^8}{{n-4} \choose {k-4}} + \frac{12n}{k^7}{{n-3} \choose {k-3}}
\end{equation}
$$
\leq n^{k-2}\left(\frac{1}{9k^8(k-4)!}+\frac{12}{k^7(k-3)!}\right)(1+o(1)).
$$
Notice that $|E_m|=(m-1)r{t \choose 2}\leq mk^4 \leq (1+o(1))\frac{n^2}{9k^8}k^4 \leq (1+o(1))\frac{n^2}{9k^4}$.
If $f \in E_m$ is independent from $e_m$ then they appear together in at most ${{n-4} \choose {k-4}}$ copies of
$K_k$. If $f$ and $e_m$ are not independent, then they may appear together in at most ${{n-3} \choose {k-3}}$ copies of
$K_k$. Suppose, w.l.o.g., that $f=(u,x)$. Since $f_{m-1}(u) \leq n/(2k)$
we know that there are at most $n/(2k)+qr(t-1) \leq n/(2k)+2k^4q$ choices for $f$ where $q$ is the number of edges
$e_j$ with $j < i$ and which have $u$ as an endpoint. However,  $q \leq \Delta(A_\pi) \leq6n/k^7$. Thus, we have that
\begin{equation}
\label{e6} 
a_4 \leq (1+o(1))\frac{n^2}{9k^4}{{n-4} \choose {k-4}}+ \left(\frac{n}{2k}+2k^4\frac{6n}{k^7}\right){{n-3} \choose {k-3}}
\end{equation}
$$
\leq n^{k-2}\left(\frac{1}{9k^4(k-4)!}+\frac{1}{2k(k-3)!}+\frac{12}{k^3(k-3)!}\right)(1+o(1)).
$$
In order to estimate $a_5$ we need to estimate the size of $V_m$. Since $|E_m| \leq (1+o(1))\frac{n^2}{9k^4}$
we have that $|V_m|(n/(2k)-r(t-1)) < (1+o(1))2n^2/9k^4$. Thus, $|V_m| \leq (1+o(1))4n/9k^3$.
Trivially, each vertex appears it at most $(n-1)/(t-1)$ elements of $F_\pi$. Thus, the overall number of elements of $F_\pi$
containing an element of $V_m$ is at most $(1+o(1))4n^2/(9k^3(t-1))$. It follows that
$|F_m| \leq (1+o(1)){t \choose 2}4n^2/(9k^3(t-1))$.
As before, for $f \in F_m$ which is independent of $e_m$ there are at most ${{n-4} \choose {k-4}}$ copies of $K_k$
containing both $f$ and $e_m$.
If $f \in F_m$ is not independent with $e_m$ then assume $f=(u,x)$. We therefore must have some $y \in V_m$
(possibly $y=x$) such that $(u,y)$ and $(u,x)$ are in the same element of $F_\pi$. It follows that there are at
most $2|V_m|(t-1)$ choices for $f$ which shares an endpoint with $e_m$ and, as before, each such $f$ appears
together with $e_m$ in at most ${{n-3} \choose {k-3}}$ copies of $K_k$. We therefore have
\begin{equation}
\label{e7}
a_5 \leq (1+o(1)){t \choose 2}\frac{4n^2}{9k^3(t-1)}{{n-4} \choose {k-4}}+2(t-1)(1+o(1))\frac{4n}{9k^3}{{n-3} \choose {k-3}}
\end{equation}
$$
\leq n^{k-2}\left( \frac{4}{9k^2(k-4)!} +\frac{16}{9k^2(k-3)!}  \right)(1+o(1)).
$$
By inequalities (\ref{e3}), (\ref{e4}),(\ref{e5}),(\ref{e6}) and (\ref{e7}) we have that $a_1 > a_2+a_3+a_4+a_5$ since
$$
0.95 > \frac{(k-2)(k-3)}{9k^8} + \frac{12(k-2)}{k^7}+\frac{(k-2)(k-3)}{9k^4}+\frac{k-2}{2k}+\frac{12(k-2)}{k^3}
+\frac{4(k-2)(k-3)}{9k^2}+\frac{16(k-2)}{9k^2}
$$
holds for all $k \geq 3$.
\endpf

\noindent
We have now completed the proof of Theorem \ref{t21}. \endpf

\section{Concluding remarks and open problems}
\begin{itemize}
\item
Although Theorem \ref{t1} and Theorem \ref{t2} are stated only for $K_k$, it is easy to see that these theorems also hold for any
$k$-vertex graph. Indeed, if $H$ has $k$ vertices then, trivially, $K_k$ has a fractional $H$-decomposition. Thus, any graph
which has a fractional $K_k$-decomposition also has a fractional $H$-decomposition. It follows that Theorem \ref{t2} holds also for
$H$. Combining Theorem \ref{t3} and Theorem \ref{t2} as explained in the introduction gives that  Theorem \ref{t1} also holds
for $H$.
\item
For all $k \geq 3$, Theorem \ref{t1} gives a lower bound of $1/9k^{10}$ for $c_{K_k}$.
We now prove that  $c_{K_k} \leq 1/(k+1)$. We will prove something slightly stronger;
For all $\epsilon > 0$ there exists $\delta > 0$ and a graph $G$ with $n$ vertices and $\delta(G) \geq n(1-1/(k+1))-\epsilon n$
for which $\nu^*_{K_k}(G) \leq (1-\delta)(e(G)/e(K_k))$. We will use a modification of a construction from \cite{Gu}
for this purpose. Let $s$ be a positive integer and let $H_s$ be any $r$-regular graph with $2s(k^3-k)$ vertices and
$r=4s(k^2-k)-d$ where $d = \lfloor \epsilon 2s(k^3-k)(k-1) \rfloor$. Such graphs clearly exist for $s$ sufficiently large
as a function of $k$ and $\epsilon$. Let $G_s$ be the graph constructed by blowing up each vertex of $K_{k-1}$ to a copy
of $H_s$. Clearly, $G_s$ has $n=2s(k^3-k)(k-1)$ vertices. $G_s$ is regular of degree $\delta=r+2s(k^3-k)(k-2)$.
Notice that $\delta = n(1-1/(k+1)) - \lfloor \epsilon n \rfloor$.
However, any $K_k$ in $G$ must contain an edge from one of the blown up copies of $H_s$. It follows that
$$
\nu^*_{K_k}(G) \leq (k-1)e(H_s)=(k-1)s(k^3-k)(4s(k^2-k)-d).
$$
But
$$
e(G)=(2s(k^3-k))^2 {{k-1} \choose 2} + (k-1)s(k^3-k)(4s(k^2-k)-d).
$$
It follows that $\nu^*_{K_k}(G) \leq (1-\delta)(e(G)/{k \choose2})$ where $\delta=\delta(\epsilon,k)$.
\item
Theorem \ref{t2} can be implemented in polynomial time, since given an input graph $G$
with $\delta(G) \geq (1-1/9k^{10})n+o(n)$ we are guaranteed by Lemma \ref{t2} that the solution
to the linear program computing $\nu^*_{K_k}(G)$ is $e(G)/e(K_k)$. It is shown in \cite{HaRo} that
Theorem \ref{t3} can be implemented in polynomial time, namely, any fractional packing can be converted in
polynomial time to an integral packing whose value differs
from $\nu^*_{K_k}(G)$ by $o(n^2)$. It follows that Theorem \ref{t1} can also be implemented in
polynomial time. The same arguments hold if we replace $K_k$ with any $k$-vertex graph $H$.
\item
It is plausible Theorem \ref{t1} can also be proved for hypergraphs.
Namely, for any positive integers $k$  and $r$ with $k>r$ there exists $\delta=\delta(k,r)$ such that the following
holds. Let $K_k^r$ denote the complete $r$-uniform hypergraph on $k$ vertices. Then, any $n$-vertex $r$-uniform hypergraph
$H$ with minimum degree at least ${{n-1} \choose {r-1}}(1-\delta(k,r))$ has a $K_k^r$-packing of size at least
$(1-o(1))e(H)/e(K_k^r)$. We already have some partial results in this direction. The proof in the hypergraph case turns out
to be significantly more involved than in the graph-theoretic case. Details will appear in a separate paper.
\item
The constant $1/9k^{10}$ is chosen to accommodate all $k \geq 3$ throughout all the lemmas. By carefully reviewing all
computations for the case $k=3$ it is easy to get $c_{K_3} \geq 1/90000$ which is 6 times better than the general constant.
We do not bother with the details since there is no indication that this improved lower bound is close to the truth.
In fact, it the conjecture of Nash-Williams \cite{Na} is true then $c_{K_3}=1/4$.
\item
Theorem \ref{t1} should be compared to other packing results which guarantee a tight packing. If all edges of the graph $G$
lie on $\alpha n^{k-2}(1+o(1))$ copies of $K_k$ for some $\alpha$ then the result of Frankl and R\"odl \cite{FrRo} guarantees a packing of size
$(1-o(1))e(G)/e(K_k)$. However, for any $\delta > 0$, there are graphs with minimum degree
at least $(1-\delta)n$ for which no such $\alpha$ exists. For any $\delta > 0$, applying the result of Kahn \cite{Ka} to graphs with
minimum degree at least $(1-\delta)n$ yields a packing of size $(1-\beta)e(G)/e(K_k)$, where $\beta$ is a {\em constant} that depends
on $\delta$ and $k$ and tends to zero as $\delta$ tends to zero.
\item
Theorem \ref{t21} gives a nontrivial minimum degree requirement which guarantees the existence of an ${\cal F}$-decomposition for the
family ${\cal F}=\{K_k,K_{2k-1},K_{2k-1}^-\}$. It is interesting to find other more general families ${\cal F}$ for which
nontrivial minimum degree conditions guarantee an ${\cal F}$-decomposition, and which do not rely on
the horrible bounds from \cite{Gu}. Notice that if  ${\cal F}$ contains a bipartite graph with minimum degree one this problem is
solved in \cite{Yu}.
\end{itemize}

\end{document}